\documentclass[a4paper,12pt]{article}
\usepackage{amsmath}
\usepackage{latexsym}
\usepackage{amsfonts}
\usepackage{amssymb}
\usepackage{color}
\usepackage[latin1]{inputenc}
\usepackage[english]{babel}
\usepackage{natbib}

\newtheorem{proposition}{Proposition}

\newenvironment{proof}[1][Proof]{\textbf{#1.} }{\ \rule{0.5em}{0.5em} \medskip \\}
\newtheorem{remark}{Remark}

\newcommand{\new}{\newcommand}

\new{\bg}{\begin}
\new{\iii}{\begin{enumerate}}
\new{\fff}{\end{enumerate}}
\new{\mfi}{\begin{eqnarray*}}
\new{\mff}{\end{eqnarray*}}
\new{\mfni}{\begin{eqnarray}}
\new{\mfnf}{\end{eqnarray}}
\new{\beeq}[2]{\begin{equation}\label{#1}{#2}\end{equation}}
\new{\eqn}[1]{~(\ref{#1})}

\new{\room}{\ \ \ \ }
\new{\tx}[1]{\mathrm{#1}}

\new{\emme}{L^0(X,\mu)}

\new{\de}{{\rm d}}

\new{\eps}{\epsilon}
\new{\ga}{\gamma}
\new{\Ga}{\Gamma}
\new{\la}{\lambda}

\new{\nor}[1]{\left\|{#1}\right\|}
\new{\norh}[1]{{\left\|{#1}\right\|_\hh}}
\new{\norq}[1]{{\left\|{#1}\right\|_q}}
\new{\norp}[1]{{\left\|{#1}\right\|_p}}
\new{\noru}[1]{{\left\|{#1}\right\|}}
\new{\normi}[1]{\left\|{#1}\right\|_{\infty}}
\new{\scal}[2]{\left\langle{#1},{#2}\right\rangle}
\new{\scalh}[2]{{\left\langle{#1},{#2}\right\rangle}_\hh}
\new{\set}[1]{\{{#1}\}}

\new{\runo}{{\mathbb R}}
\new{\cuno}{{\mathbb C}}
\new{\nat}{{\mathbb N}}

\new{\hh}{{\mathcal H}}
\new{\hhg}{{\hh_\ga}}
\new{\hhm}{{\hh_\mu}}
\new{\kk}{\mathcal K}

\new{\bo}[1]{{\mathcal B}(#1)}
\new{\cc}[1]{{\mathcal C}(#1)}
\new{\ccc}[1]{{\mathcal C}_c(#1)}

\new{\luno}{L^1(X,\mu)}
\new{\ldue}{L^2(X,\mu)}
\new{\lpi}{L^p(X,\mu)}
\new{\lqu}{L^q(X,\mu)}
\new{\lri}{L^r(Y,\nu)}
\new{\lduet}{L^2(Y,\nu)}

\new{\maaa}{\mu\mathrm{-a.a.}}
\new{\lunot}{L^1(X)}
\new{\lduenot}{L^2(Y)}
\new{\ldueXnot}{L^2(X)}

\new{\supp}[1]{\operatornamewithlimits{supp}\,#1}
\new{\lin}[1]{\operatornamewithlimits{span}\{\,#1\}}
\new{\range}[1]{\operatornamewithlimits{Im}\,#1}
\new{\Ker}[1]{\operatornamewithlimits{ker}\,#1}

\new{\mae}{\mu\mathrm{-a.e.}}
\new{\R}{R}

\begin{document}

\author{C.~Carmeli\thanks{C.~Carmeli,
Dipartimento di Fisica, Universit\`a di Genova, and I.N.F.N.,
Sezione di Genova, Via Dodecaneso~33, 16146 Genova, Italy. e-mail:
carmeli@ge.infn.it}, E. De Vito\thanks{E.~De Vito,
Dipartimento di Matematica, Universit\`a di Modena e Reggio
Emilia, Via Campi 213/B, 41100 Modena, Italy, and I.N.F.N., Sezione di
Genova, Via Dodecaneso~33, 16146 Genova, Italy. e-mail: devito@unimo.it}, A.
Toigo\thanks{A.~Toigo,
Dipartimento di Fisica, Universit\`a di Genova, and I.N.F.N.,
Sezione di Genova, Via Dodecaneso~33, 16146 Genova, Italy. e-mail:
toigo@ge.infn.it},
}
\title{Reproducing kernel Hilbert spaces and Mercer theorem}
\date{\today}

\maketitle

\begin{abstract}
We characterize the reproducing kernel Hilbert spaces whose
elements are $p$-integrable functions in terms of the {\em
boundedness} of the integral operator whose kernel is the
reproducing kernel. Moreover, for $p=2$ we show that the
spectral decomposition of this integral operator gives a complete
description of the reproducing kernel.
\end{abstract}
\section{Introduction}

In recent years there is a new interest for the theory of reproducing
kernel Hilbert spaces in different frameworks, like statistical learning
theory \citep{cusm}, signal analysis \citep{dau} and quantum
mechanics \citep{ali}. In particular, for these applications there is
often the need of reproducing kernel Hilbert spaces having some additional regularity
property, as continuity or square-integrability.\\
This paper is both a research article and a self-contained survey about the characterization of
the reproducing kernel Hilbert spaces whose elements are  continuous, measurable or
$p$-integrable functions ($1\leq p\leq \infty$). As
briefly reviewed in Section~\ref{pre}, this problem is
equivalent to study the weak regularity properties of maps
taking values in an arbitrary Hilbert space \citep{saitoh88}. In the following sections, we take this last
point of view  and, given a set $X$ and a Hilbert space $\hh$, 
we characterize the weak regularity properties of a map $\ga:X\to\hh$ in terms of the corresponding properties of the
associated kernel
\beeq{a}{X\times X \ni (x,y)\mapsto \scalh{\ga_y}{\ga_x}\in \cuno.}
More precisely, in Section~\ref{sec_meas} we prove that, if $X$ is a
measurable set and $\hh$ is a separable Hilbert space, $\ga$ is weakly
measurable if and only if the associated kernel\eqn{a} is separately
measurable. Our proof is
an easy consequence of the equivalence between weak and strong
measurability. Moreover, if $\hh$ is a space of square-integrable
functions, we recall a result of \citet{pettis} clarifying the
relation between vector valued maps and the theory of integral
operators \citep{hasu}.\\
In Section~\ref{sec_int} we show that, if $X$ is a measurable set endowed with a $\sigma$-finite
measure $\mu$ and $\hh$ is a separable Hilbert space,  $\ga$ is weakly $p$-integrable if and only if the
integral operator with kernel given by Eq.\eqn{a} is bounded from $L^{\frac{p}{p-1}}(X,\mu)$ into
$L^p(X,\mu)$ where $1\leq p\leq\infty$. Up to our knowledge this result is new, though is
deeply based on the theory of Pettis integral \citep{hiph}. Moreover, for $p=1$ we prove that the operator 
\[\hh \ni v\mapsto \scalh{v}{\ga(\cdot)}\in L^1(X,\mu)\]
is always compact. For finite measures this result is due to
\citet{pettis}, but it appears original for non finite measures. A
brief discussion of the compactness for $1<p<\infty$ is also given to show
that the strong $p$-integrability of the map $\ga$ is a sufficient condition,
but is not necessary.\\ 
In Section~\ref{sec_cont}, if $X$ is a locally compact
space, it is proved that $\ga$ is weakly continuous if and only if the kernel\eqn{a} is locally
bounded and separately continuous, whereas the continuity of the
kernel on $X\times X$ is equivalent to the compactness of the operator  
\[\hh \ni v\mapsto\scalh{v}{\ga(\cdot)}\in \cc{X}.\]
The results we present are due to \citet{schwartz} in the
framework of reproducing kernel Hilbert spaces and we review them
giving an elementary proof based on standard functional analysis
tools.\\ 
Finally, in Section~\ref{sec_mercer} if $X$ is a locally compact second countable Hausdorff space endowed with a
positive Radon measure $\mu$ and $\hh$ is a
reproducing kernel Hilbert space such that $\hh\subset
L^2(X,\mu)\cap\cc{X}$, we characterize the space $\hh$ and the
reproducing kernel $\Ga$ in terms of the spectral decomposition of the
integral operator of kernel $\Ga$. When $X$ is compact, this kind of
result is known as {\em Mercer theorem} \citep{hoch}. Extensions of Mercer theorem can be found
in \citet{novi} and references therein. Our proof is very simple and
general since it is based on the polar decomposition of the canonical
inclusion of $\hh$ into $\ldue$.

\section{Notations}\label{pre}
In this section we fix the notation, give the main definitions and
review the connection between reproducing kernel Hilbert spaces and
vector valued maps.  

If $E$ is a (complex) Banach
space, $\nor{\cdot}_E$ denotes the norm of $E$, and $E^*$ is
the Banach space of continuous antilinear functionals of $E$.
We let $\scal{\cdot}{\cdot} : E^*\times E\to \cuno$
be the canonical pairing.
If $F$ is another Banach space and $A:E\to
F$ is a bounded linear operator, then $A^*:F^*\to
E^*$ is the adjoint of $A$. If
$\hh$ is a Hilbert space, $\scalh{\cdot}{\cdot}$ denotes the scalar
product, linear in the first argument and, by means of the scalar product, $\hh^*$
is canonically identified with $\hh$.

If $X$ is a set, $\cuno^X$ is the vector space of  all the complex
functions on $X$. \\
If $X$ is a measurable space endowed with a $\sigma$-finite positive measure
$\mu$, $\emme$ denotes the topological vector space
of all  measurable  complex functions\footnote{As usual, a function is
  identified with its equivalence class $\mu$-almost everywhere.} on $X$
endowed with the topology of the convergence in measure on subsets of
finite measure \citep{schwartz1}. Given $1\leq p <\infty$, $\lpi$ is the 
Banach space of functions  $f\in\emme$ such that $\|f\|^p$
is $\mu$-integrable, and $L^\infty(X,\mu)$ is
the Banach space of elements $f\in\emme$ that are bounded $\mu$-almost everywhere. 

If $X$ is a locally compact Hausdorff space, $\cc{X}$ denotes the
space of continuous functions on $X$ endowed with the open-compact
topology \citep{kelley}. If $X$ is second countable locally compact Hausdorff space,
a {\em positive Radon measure} on $X$ is a positive measure $\mu$
defined on the Borel $\sigma$-algebra of $X$ and finite on compact
subsets, and $\supp{\mu}$ denotes the support of $\mu$.

Given a map $\ga$ from a set $X$ into a Hilbert space 
$\hh$, we denote by 
$A_\ga: \hh \to \cuno^X$ the linear operator 
\beeq{acar}{(A_\ga v)(x)=\scalh{v}{\ga_x}\room \forall x\in X\  v\in\hh,}
by $\Ga: X\times X\to\cuno$ the kernel
\beeq{Gamma}{\Ga(x,t)=\scalh{\ga_t}{\ga_x}=A_\ga(\ga_t)(x)\room
  \forall x,t\in X}
and we let
\beeq{hhg}{\hhg=\overline{\rm span}\,\{\ga_x\in\hh\,\vert\,
      x\in X\},}
where $\overline{\rm span}$ is the closure of the linear span.

We now recall some basic definitions about {\em weak regularity} properties of vector
valued maps.
\bg{definition}\label{def_car} 
Let $X$ be a set, $\hh$ a Hilbert space and $\ga:X\to\hh$.  
\iii
\item Assume $X$ be a measurable space. The map $\ga$ is  {\em
    weakly measurable} if the function $A_\ga v$ is measurable for all $v\in\hh$.
\item Assume $X$ be a measurable space endowed with a $\sigma$-finite
  measure $\mu$ and $1\leq p\leq \infty$. The map $\ga$ is  {\em
  weakly $p$-integrable} if the function $A_\ga v$ is $p$-integrable for all 
$v\in\hh$. 
\item Assume that $X$ is a locally compact Hausdorff space. The map
  $\ga$ is {\em weakly continuous} if the function $A_\ga v$ is
  continuous for all $v\in\hh$.
\fff
\end{definition}
The assumptions that $\mu$  is
a $\sigma$-finite measure and $X$ is locally compact will avoid technical problems.
Moreover, the above properties clearly depend only on $\hhg$.  
If $X$ is a measurable space with a $\sigma$-finite measure $\mu$ and
$\ga$ is measurable, we let 
\beeq{essrange}{\mathcal S=\set{v\in\hh \vert \mu(\ga^{-1}(B(v,\eps)))>0 \ \forall
  \eps>0}\room\room \hhm=\overline{\rm span}\,\mathcal S\ ,}
 where $B(v,\eps)=\set{w\in\hh \,\vert\, \norh{w-v}<\eps}$. The closed
 subset $\mathcal S$ is the {\em essential range} of $\ga$
 and $\mathcal S\subset\hhg$ so that $\hhm\subset\hhg$. 

We need also the definition of kernel of positive type.
\bg{definition}\label{positive_type}
Given a set $X$, a complex kernel $\Ga:X\times X\to\cuno$ is called
of positive type if
\beeq{positive}{\sum_{i,j}^\ell c_i \overline{c_j}\Ga(x_i,x_j) \geq 0}
for any $\ell\in\nat$, $x_1,\ldots,x_\ell\in X$ and $c_1,\ldots,c_\ell\in\cuno$.
\end{definition}
\begin{remark}
In the complex case, the positivity condition (\ref{positive})
ensures that
\beeq{symmetric}{\Ga(x,t)=\overline{\Ga(t,x)}\room\forall x,t\in X.}
This is no longer true in the real case. In this case, a kernel
$\Ga:X\times X\to\mathbb{R}$ is called of positive type if
\mfi 
\Ga(x,t) & = & \Ga(t,x)\room \forall x,t\in X \\
\sum_{i,j}^\ell c_i c_j \Ga(x_i,x_j) & \geq & 0 
\room \forall\ell\in\nat ,\ x_1,\ldots,x_\ell\in X,\ c_1,\ldots,c_\ell\in\mathbb{R}.
\mff
\end{remark}
In the framework of harmonic analysis, kernel of positive type are
also called positive definite, and, in the context of reproducing
kernel Hilbert spaces, Aronszajn kernel.

We now recall the definition of {\em bounded kernel} from the
theory of integral operators.
\bg{definition}\label{bound}
Let $X$ and $Y$ two measurable spaces endowed with $\sigma$-finite measures $\mu$ and
$\nu$, respectively. Fixed $1\leq p,r\leq \infty$, a measurable kernel
$K:X\times Y\to\cuno$ is called  {\em $(r,p)$-bounded} if, for
all $\phi\in \lri$, 
\iii
\item there is a $\mu$-null set $X_\phi\subset X$
such that  the function  $K(x,\cdot)\phi(\cdot)$ is in
$L^1(Y,\nu)$ for all $x\not\in X_\phi$; 
\item the map 
$$
x\mapsto \int_Y K(x,y)\phi(y)\,\de\nu(y)
$$ 
is in $\lpi$.
\fff
Moreover, if the following condition holds
\iii
\item[1')] for $\mu$-almost all $x\in X$, $K(x,\cdot)\in
  L(Y,\nu)^{\frac{r}{r-1}}$ (if $r=1$, $\frac{r}{r-1}=\infty$ and, if
  $r=\infty$, $\frac{r}{r-1}=1$), 
\fff
$K$ is called a {\em Carleman $(r,p)$-bounded} kernel. 
\end{definition}
Since $\nu$ is $\sigma$-finite, Condition~1') is
equivalent to Condition~1) and the fact that the $\mu$-null set
$X_\phi$ can be chosen in such a way to be independent of $\phi$. 
In \citet{hasu} there is an example of $(2,2)$-bounded
kernel, which is not a Carleman kernel.
In the above definition, {\em boundedness} refers to the fact that the
integral operator of kernel $K$ is bounded from $\lri$ to $\lpi$, as shown in 
Proposition~\ref{bounded}.

Finally, we recall that a reproducing kernel Hilbert space $\hh$ on $X$ 
is a subspace of $\cuno^X$ such that $\hh$ is a Hilbert space and, for all $x\in X$, 
there is
a function $\ga_x\in\hh$ satisfying
\beeq{valuta}{ f(x) =  \scalh{f}{\ga_x} \room \forall f\in\hh.}  
The corresponding reproducing kernel of $\hh$ is defined by
\beeq{reprokernel}{\Ga(x,t) = \scalh{\ga_t}{\ga_x} = \ga_t (x)\room
  x,t\in X.}  
We now review the connection between reproducing kernel Hilbert spaces, Hilbert space valued maps
and kernels of positive type. The following result has been obtained by
many authors, see
\citet{aron,godement,kolmogorov,krein1,krein2,schoenberg2} and, for a complete list of references,
\citet{bekka,hille,saitoh88,saitoh97,schwartz}.
\bg{proposition}\label{range}
Let $X$ be a set.
\iii
\item Given a kernel $\Ga:X\times X\to\cuno$ of
positive type, there is a unique reproducing kernel Hilbert space
$\hh$ on $X$ with reproducing kernel $\Ga$. The inclusion of $\hh$ into
$\cuno^X$ is the operator $A_\ga$ associated with the map $\ga:X\to
\hh$ defined by
\beeq{ga_rkhs}{\ga_x=\Ga(\cdot,x)\room x\in X.}
The kernel\eqn{Gamma} associated with $\ga$ is precisely the
reproducing kernel of $\hh$ and
\beeq{genera}{\hh =\hhg\ .}
\item Given a Hilbert space $\hh$ and a map $\ga:X\to\hh$,  
the associated kernel $\Ga$ is of positive type 
and $\ker{A_\ga}=\hhg^\perp$. In particular, $A_\ga$ is a unitary
operator from $\hhg$ onto the reproducing kernel 
Hilbert space with reproducing kernel $\Ga$.
\fff
\end{proposition}
\bg{proof}
We report the proof given in \citet{bekka}.
\iii
\item Let $\ga_x=\Ga(\cdot,x)\in\cuno^X$ with $x\in X$ and
$$\hh_0=\lin{\ga_x \,\vert\, x\in X}\subset\cuno^X.$$
If $f=\sum_i c_i \ga_{x_i}\in\hh_0$ and  $g=\sum_j
d_j\ga_{t_j}\in\hh_0$, the definition of $\ga_x$ implies that
$$\sum_{ij} c_i \overline{d_j} \Ga(t_j,x_i)  =  \sum_{j}
\overline{d_j} f(t_j) = \sum_{i} c_i \overline{g(x_i)},$$
so the following sequilinear form on $\hh_0\times\hh_0$
$$(f,g)\mapsto \scal{f}{g}:= \sum_{ij} c_i \overline{d_j} \Ga(t_j,x_i) $$ 
is well defined. Equation\eqn{symmetric} ensures that
$\scal{\cdot}{\cdot}$ is hermitian and Eq.\eqn{positive} that
$\scal{f}{f}\geq 0$ for all $f\in\hh_0$.
Let now $x\in X$, the choice $g=\ga_x$ in the above equation implies that
\beeq{w1}{f(x)=\scal{f}{\ga_x}\room\forall x\in X,}
for all $f\in\hh_0$. Assume now that $\scal{f}{f}=0$, the
Cauchy-Schwarz inequality gives
$$ |f(x)|\leq \sqrt{\scal{f}{f}}\sqrt{\scal{\ga_x}{\ga_x}}=0, $$
and, hence, $f=0$. This shows that $\scal{\cdot}{\cdot}$ is a scalar
product on $\hh_0$. Finally, if $(f_n)_{n\in\nat}$ is a Cauchy sequence in $\hh_0$,
Eq.\eqn{w1} implies that, for all $x\in X$, $(f_n(x))_{n\in\nat}$
converges to $f(x)\in\cuno$, so that the completion of $\hh_0$ is a
subspace of $\cuno^X$ and Eq.\eqn{w1} holds for all $f\in\hh$, showing
that $\hh$ is a reproducing kernel Hilbert space with kernel
$\Ga$. The uniqueness is evident.  Finally, let $\ga$ as in
Eq.\eqn{ga_rkhs}, Eq.\eqn{valuta} implies that the corresponding
operator $A_\ga$ is the inclusion of $\hh$ into $\cuno^X$. 
Equation\eqn{genera} follows by density of $\hh_0$.
\item A simple check shows that the kernel $\Ga$ associated with $\ga$
  is of positive type and that
  $\ker{A_\ga}=\hhg^\perp$. In particular, $A_\ga$ is a bijective map
  from $\hhg$ onto $\range{A_\ga}$. Let now $\hh_\Ga$ be the reproducing kernel 
Hilbert space 
  with kernel $\Ga$ given by point 1 of this proposition. Both $\range{A_\ga}$ and 
$\hh_\Ga$ are subspaces of
$\cuno^X$. In particular, for all $x\in X$, the function
  $\Ga(\cdot,x)=A_\ga \ga_x$ is both in $\range{A_\Ga}$ and in
$\hh_\Ga$, see Eq.\eqn{acar} and Eq.\eqn{reprokernel},  respectively.  
Let now $x,t\in X$, the cited equations give
$$\scal{A_\ga\ga_t}{A_\ga\ga_x}_{\hh_\Ga}=\scal{\Ga(\cdot,t)}{\Ga(\cdot,x)}_{\hh_\Ga}=
\Ga(x,t)=\scalh{\ga_t}{\ga_x},$$
showing that $A_\ga$ is an isometry from $\lin{\ga_x\,|\,x\in X}$
onto $\lin{\Ga(\cdot, x)\,|\,x\in X}$. The claim follows
by definition of $\hhg$ and Eq.\eqn{genera} applied to $\hh_\Ga$.
\fff
\end{proof}
The construction of the Hilbert space $\hh$ given the kernel $\Ga$ is
related to the GNS construction and it is also known as Kolmogorov
theorem \citep{kolmogorov}. In the above proof, $\hh$ is directly defined
as a subspace of $\cuno^X$. In doing so, one has to prove that the
scalar product is well defined, but, due to Eq.\eqn{valuta}, it is
strictly positive, so that there is no need to quotient with respect to the
vectors of null norm. Another possibility is to define $\hh$ as a
subspace of the formal linear combinations of elements
$\Ga(\cdot,x_i)$, see, for example, \citet{dutkay}.

Proposition~\ref{range} shows that there is a one-to-one correspondence
between Hilbert space valued maps, complex kernels of positive type
and reproducing kernel Hilbert spaces. In particular, given a map
$\ga$, the Hilbert space where the map takes value can be identified
by means of the operator $A_\ga$ with a unique reproducing kernel
Hilbert space, which is a subspace of $\cuno^X$. Conversely, any
reproducing kernel Hilbert space defines uniquely a Hilbert space
valued map $\ga$ such that the operator $A_\ga$ is the inclusion into
$\cuno^X$. In both cases, the weak regularity properties of $\ga$
stated in Definition~\ref{def_car} are in correspondence with the regularity properties of the 
complex functions
in the reproducing kernel Hilbert space. In the following, we will show that
these properties can be completely characterized in terms of the kernel $\Ga$. 
Due to the above equivalence, we will state the results for an
arbitrary Hilbert space valued map.
Here, for the convenience of the reader,  we
summarize them in the case of a separable reproducing kernel Hilbert space $\hh$
and for $X$ being a locally compact second countable Hausdorff space endowed with a
positive Radon measure $\mu$.
\iii
\item  Proposition~\ref{misurabilita}: $\hh\subset\emme$ if and only if $\Ga$ is 
separately measurable;
\item  Proposition~\ref{p-integrabilita}: $\hh\subset\lpi$ if and only if $\Ga$ is a
  $(\frac{p}{p-1},p)$-bounded kernel, for all $1\leq p\leq \infty$;
\item Proposition~\ref{continuita}: $\hh\subset\cc{X}$ if and only if $\Ga$ is 
separately continuous
  and locally bounded.
\fff
The assumption that $\hh$ is separable is essential (see
Example~\ref{nonseparabile}) and, in general, the inclusion of $\hh$
into $\emme$ is not injective. However, if $\ga$ is weakly continuous and $\supp{\mu}=X$, $\hh$ is
separable and the inclusion is injective (see
 Eqs.\eqn{genera},\eqn{kernel},\eqn{kernel-bello}).

\section{Measurability}\label{sec_meas}
In this section we characterize the weak measurable maps and we
discuss the relation with the theory of integral operators. The following proposition  
is based on the well-known equivalence between weak
and strong measurability for maps with separable range (see, for example, \citet{hiph,pettis}).
\bg{proposition}\label{misurabilita}
 Let $X$ be a measurable space, $\hh$ a Hilbert space and $\ga:X\to\hh$. 
Assume that $\hhg$ is separable, then the following conditions are equivalent:
\iii
\item the map $\ga$ is weakly measurable;
\item the map $\ga$ is measurable from $X$ to $\hh$;
\item the function $\Gamma$ is measurable from $X\times X$ 
into $\cuno$;
\item for all $x\in X$, the  function $\Ga(x,\cdot)$ is measurable
  from $X$ into $\cuno$.
\fff
 If $X$ is endowed with a
$\sigma$-finite measure $\mu$ and
one of the above conditions holds, then $\range{A_\ga}\subset\emme$,
  the operator $A_\ga$ is a continuous from $\hh$ into $\emme$ and
\beeq{kernel}{\Ker{A_\ga}=\hhm^\perp\ .}
\end{proposition}
\bg{proof}{~}
\iii
\item[$1)\Rightarrow 2)$] Since $\hhg$ is separable, there is a denumerable
Hilbert basis $(e_n)_{n\in I}$ of $\hhg$ and
$$ \ga_x = \sum_{n\in I} \overline{(A_{\gamma}e_n)(x)}\,
e_n\room\forall x\in X.$$
By definition $A_\ga e_n=\scal{e_n}{\ga(\cdot)}$ and is
measurable by assumption and $\ga$ is measurable since $I$ is
denumerable. 
\item[$2)\Rightarrow 3)$] Since scalar product is continuous and $\hhg$ is separable,
the function $\Gamma(x,t) =\scalh{\ga_t}{\ga_x}$ is measurable.
\item[$3)\Rightarrow 4)$] Given $x\in X$, the map $t\mapsto(x,t)$ is
  clearly measurable, so Condition~(4) follows.
\item[$4)\Rightarrow 1)$] Let $v=\sum_{i=1}^N c_i \ga_{x_i}$ with
  $x_i\in X$. By assumption $\Ga(\cdot,x_i)$ is measurable for all
  $i$, so is  $A_\ga v = \sum_i c_i \Ga(\cdot,x_i)$. By density and Eq.\eqn{hhg},
  $A_\ga v$ is measurable for all $v\in \hhg$. If $v\in\hhg^\perp$, $A_\ga
  v=0$ and, hence, is measurable.
\fff
Assume now one of the four equivalent conditions. Since
pointwise convergence implies convergence in measure on subsets of
finite measure \citep{schwartz1}, Eq.\eqn{acar} gives that $A_\ga$
is continuous from $\hh$ to $\emme$. 

We now prove $\Ker{A_\ga}=\hhm^\perp$. Since
$\hhg^\perp\subset\Ker{A_\ga}$, without loss of generality, we can assume
that $\hhg=\hh$ with $\hh$ separable and, by definition of $\hhm$, it is
enough to show that $\Ker{A_\ga}=\mathcal S^\perp$, where $\mathcal S$
is the essential range of $\ga$.\\
The set $X_0=\set{x\in X\,|\, \ga_x\not\in \mathcal S}$
is of null measure. Indeed, since $\mathcal S$ is closed and $\ga$ is (strongly)
measurable, $X_0$ is measurable. Moreover,  by definition, if $x\in X_0$, there is
$\eps_x >0$ such that $\mu(\ga^{-1}(B(\ga_x,\eps_x)))=0$.
Since $\hhg$ is separable, each member of the family
$\set{B(\ga_x,\eps_x)}_{x\in X_0}$ is a union of members of a fixed countable
family of open sets $\set{\Omega_i}_{i\in\nat}$, with each $\Omega_i$
contained in $B(\ga_x,\eps_x)$ for some $x\in X_0$. Clearly
$\mu(\ga^{-1}(\Omega_i))=0$ and $X_0\subset \cup_{i\in\nat}
\ga^{-1}(\Omega_i)$, so that $\mu(X_0)=0$.\\ 
Let now $v\in\mathcal S^\perp$, then $(A_\ga
v)(x)=0$ for all $x\not\in X_0$ since $\ga_x\in\mathcal S$, but
$\mu(X_0)=0$, so that $v\in\Ker{A_\ga}$. 
Conversely, let $v\in\Ker{A_\ga}$. By contradiction, assume there is
$w_0\in\mathcal S$ such that
$\scal{v}{w_0}\neq 0$. Since the scalar product is continuous, 
$$ \scal{v}{w}\neq 0\room \forall w\in B(w_0,\eps),$$  
for some $\eps>0$. In particular, $\scal{v}{\ga_x}\neq 0$ for all $x\in 
\ga^{-1}(B(w_0,\eps))$. 
However,  $\mu(\ga^{-1}(B(w_0,\eps)))>0$, since $w_0\in\mathcal S$,
so that $A_\ga v\neq 0$.  Since this last fact contradicts the
assumption, it follows that $\scal{v}{w_0}=0$ for all $w_0\in\mathcal
S$, that is, $v\in\mathcal S^\perp$.
\end{proof} 
The following example (see \citet{bourbaki}) shows that the
separability of $\hhg$ is essential in the above proposition.
\bg{example}\label{nonseparabile}
Let $X=[0,1]$ with the Lebesgue measure and $\hh$ be the Hilbert space of
functions $v:[0,1]\to\cuno$ such that
$$\sum_{x\in X} |v(x)|^2 <+\infty,$$
with scalar product 
$$ \scalh{v}{w} = \sum_{x\in X} v(x) \overline{w(x)}.$$
Let $A$ be a non measurable subset of $[0,1]$ and $\ga$
$$\ga_x = \left\{
\begin{array}{cc}
\delta_x & x \in A \\
0 & x\not\in A
\end{array}\right. $$
where $\delta_x(y)=0$ if $y\neq x$ and $\delta_x(x)=1$. Since
$(A_\ga v)(x)=0$ for all but denumerable $x$, $A_\ga v$ is 
measurable. However, the function
$$\Ga(x,x)= \left\{
\begin{array}{cc}
1 & x \in A \\
0 & x\not\in A
\end{array}\right. $$
is not measurable and, hence, Condition~(2) and~Condition~(3) of
Proposition~\ref{misurabilita} can not be true.
\end{example} 
We now clarify the relation with the theory
of integral operators (for a complete account see \citet{hasu}, where
weakly integrable maps are called {\em Carleman functions}). 
The following proposition, due to \citet{pettis}, characterizes the measurable maps
taking value in $L^2(Y,\nu)$ (see also Theorem 11.5 of \citet{hasu}).
For a converse result, see Proposition~\ref{op-int} below. 
\begin{proposition}\label{op-misur}
Let $X$ and $Y$ two measurable sets endowed with two $\sigma$-finite
measures $\mu$ and $\nu$, respectively.  Suppose $\lduet$ is separable.
Given a weakly measurable map
$\ga:X\to\lduet$, there exists a measurable
function $K_\ga:X\times Y\longmapsto \mathbb{C}$ and a
$\mu $-null set $N\subset X$ such that,
for $x\in X\setminus N$, one has $\gamma _{x}\left( y\right) =\overline{K_\ga\left(
x,y\right)}$ for $\nu $-almost all $y\in Y$.

The operator $A_\ga:\lduet\to\emme$ associated with $\ga$ is the integral operator
with kernel $K_\ga$ and  the kernel associated with $\ga$ is 
$$\Ga(x,x') =  \int_Y \overline{K_\ga(x',y)} K_\ga(x,y)\,\de\nu(y),$$
for $\mu\otimes\mu$-almost all $(x,x')\in X\times X$.
\end{proposition}
\begin{proof}
Let $\left( C_{n}\right) _{n\geq 1}$ be an increasing sequence of 
measurable subsets of $X$  such that $\mu(C_n)<+\infty$ and $X=\bigcup_{n}C_{n}$.
For all $n\in\mathbb N$, define by induction
\begin{eqnarray*}
A_0&=&\emptyset \\
A_{n}&=&\left\{ x\in X\mid x\in C_n,x\not\in A_{n-1}\, \tx{and}\,
\nor{\gamma_x}_{\lduenot}\leq n\text{,}\right\}
\end{eqnarray*}
Clearly, each $A_{n}$ is a subset of finite measure, $A_n\cap
A_m=\emptyset$ if $n\neq m$, and $\bigcup_{n}A_{n}=X$. Define
\begin{equation*}
\gamma _{x}^{n}=\left\{
\begin{array}{cc}
\gamma _{x} & \text{if }x\in A_{n} \\
0 & \text{elsewhere}
\end{array}
\right. \text{.}
\end{equation*}
 If $f\in L^{2}\left( X\times Y,\mu \otimes \nu \right)$, by Fubini theorem
 $f(x,\cdot)\in \lduet$ for $\mu$-almost all $x\in X$. Moreover, by separability of
$\lduet$ the map $x\mapsto \left\langle f\left( x,\cdot \right)
,\gamma _{x}^{n}\right\rangle_{\lduenot}$ is measurable.
Let $\lambda_n$ be the linear form on  $L^{2}\left( X\times Y,\mu
 \otimes \nu \right) $ given by
\begin{equation*}
\lambda _{n}\left( f\right) =\int_{X}\left\langle f\left( x,\cdot \right)
,\gamma _{x}^{n}\right\rangle _{L^{2}\left( Y\right) }\text{d}\mu \left(
x\right) \text{.}
\end{equation*}
The linear form $\lambda _{n}$ is bounded since
\begin{eqnarray*}
\left| \lambda _{n}\left( f\right) \right|  &\leq &\int_{X}\left|
\left\langle f\left( x,\cdot \right) ,\gamma _{x}^{n}\right\rangle
_{L^{2}\left( Y\right) }\right| \text{d}\mu \left( x\right)  \\
&\leq &\int_{X}\left\| f\left( x,\cdot \right) \right\| _{L^{2}\left(
Y\right) }\left\| \gamma _{x}^{n}\right\| _{L^{2}\left( Y\right) }\text{d}
\mu \left( x\right)  \\
&\leq &\left[ \int_{X}\left\| f\left( x,\cdot \right) \right\| _{L^{2}\left(
Y\right) }^{2}\text{d}\mu \left( x\right) \right] ^{1/2}\left[
\int_{X}\left\| \gamma _{x}^{n}\right\| _{L^{2}\left( Y\right) }^{2}\text{d}
\mu \left( x\right) \right] ^{1/2} \\
&\leq &n \sqrt{\mu\left( A_{n}\right)} \left\| f\right\| _{L^{2}\left( X\times
Y\right) }\text{.}
\end{eqnarray*}
By Riesz lemma, there exists unique $\Lambda _{n}\in
L^{2}\left( X\times Y,\mu\otimes\nu\right) $ such that
\begin{equation*}
\lambda _{n}\left( f\right) =\int_{X\times Y}f\left( x,y\right) \Lambda
_{n}\left( x,y\right) \text{d}\mu \left( x\right) \text{d}\nu \left(
y\right) \text{.}
\end{equation*}
 It is not restrictive to assume that
$\Lambda _{n}\left( x,\cdot \right)
\in L^{2}\left( Y,\nu \right)$ for \emph{all} $x$. Let $u\in L^{2}\left( Y,\nu\right) 
$.
For every measurable subset $C$ of $X$ with $\mu(C)<+\infty$, choosing $f\left( 
x,y\right) =\chi
_{C}\left( x\right) u\left( y\right) $, the above relations and Fubini
theorem give
\begin{eqnarray*}
\int_{C}\left\langle u,\gamma _{x}^{n}\right\rangle _{L^{2}\left( Y\right) }
\text{d}\mu \left( x\right) &=& \int_{C}\text{d}\mu \left( x\right)
\int_{Y}u\left( y\right) \Lambda _{n}\left( x,y\right) \text{d}\nu \left(
y\right)  \\
& = &\int_{C} \scal{u}{\overline{\Lambda_n(x,\cdot)}}_{\lduenot}\text{d}\mu\left( 
x\right)
 \text{,}
\end{eqnarray*}
which implies $\left\langle u,\gamma _{x}^{n}\right\rangle _{L^{2}\left(
Y\right) }= \scal{u}{\overline{\Lambda_n(x,\cdot)}}_{L^2 (Y)}$ for $\mu$-almost all 
$x$. Letting $u$ vary over a denumerable
dense subset of $\lduet $, we find that there exists a $\mu$-null set  
$N_{n}\subset X$ such that, for $x\in X\setminus N_{n}$,
\beeq{uguale}{\gamma _{x}^{n}\left( y\right) =\overline{\Lambda _{n}\left( 
x,y\right)}}
for $\nu $-almost all $y\in Y$. Hence the map
\begin{equation*}
K_\ga\left( x,y\right) =\sum_{n}\chi _{A_{n}}\left( x\right) \Lambda _{n}\left(
x,y\right) 
\end{equation*}
is finite, measurable and, for all $x\in A_n\backslash N_n$, 
$$\overline{K_\ga\left( x,y\right)} =\gamma _{x}^{n}(y)= \gamma_{x}(y)$$
for $\nu $-almost all $y\in Y$. Since $X=\bigcup_{n}A_{n}$ and
$N=\bigcup_{n}N_{n}$ is a $\mu$-null set, the first part of the
proposition follows.\\
The definition of $A_\ga$ implies that $A_\ga$ is the integral
operator of kernel $K_\ga$ and the last equation easily follows. 
\end{proof} 

\section{Integrability}\label{sec_int}
In this section we characterize weakly $p$-integrable maps. First of all, we recall
that the integral operator with a $(r,p)$-bounded kernel is continuous.  
\bg{proposition}\label{bounded}
With the assumptions of Definition~\ref{bound}, let $K$ be a
$(r,p)$-bounded kernel, then the operator $L_K:\lri\to\lpi$  
\beeq{integrale}{(L_K\phi)(x)=\int_Y K(x,y) \phi(y)\,\de\nu(y)\room \text{for}\ \maaa \ x\in X}
is bounded.
\end{proposition}
\bg{proof}
We report the proof of \citet{hasu}. Since $L_K$ is defined on all the
space $\lri$, it is enough to prove that it is closed. Let $(\phi_n)_{n\in\nat}$ be
a sequence such that it converges to $\phi$ in $\lri$ and $\psi_n=L_K\phi_n$
converges to $\psi$ in $\lpi$. Possibly passing to a double subsequence, we
can assume that both $(\phi_n)_{n\in\nat}$ and $(\psi_n)_{n\in\nat}$ converges almost 
everywhere
and that there is $\omega\in\lri$ for which $|\phi_n(y)|\leq
\omega(y)$ for $\nu$-almost all $y\in Y$ (if $r<\infty$ it is a
consequence of Fischer-Riesz theorem, see Lemma~3.9 of \citet{hasu}, for
$r=+\infty$, it is trivial). Condition~1 of Definition~\ref{bound} and the
fact that denumerable union of null sets is a null set imply that,
for $\mu$-almost all $x\in X$, the sequence 
$\left(K(x,\cdot)\phi_n(\cdot)\right)_{n\in\nat}$
converges to $K(x,\cdot)\phi(\cdot)$ for $\nu$-almost all $y\in Y$ and it is
bounded almost everywhere by $|K(x,\cdot)\omega(\cdot)|\in
L^1(Y,\nu)$. The dominated convergence theorem gives that
$$\lim_{n\to+\infty}\psi_n(x)=\lim_{n\to+\infty}\int_Y K(x,y)\phi_n(y)\,\de\nu(y)= 
\int_Y
K(x,y)\phi(y)\,d\nu(y)$$
for $\mu$-almost all $x\in X$, so that $\psi(x)= \int_Y K(x,y)\phi(y)\,d\nu(y)$ in 
$\lpi$.   
\end{proof}

We are now in position to state the main result of the paper.
\bg{proposition}\label{p-integrabilita}
Let $X$ be a measurable space, $\mu$ a $\sigma$-finite measure on $X$ and $\hhg$ a 
separable Hilbert space.
Let $1\leq p \leq\infty$, then the following conditions are equivalent 
\iii
\item  the map $\ga$ is weakly $p$-integrable;
\item  the kernel $\Ga$ is $(q,p)$-bounded with $q=\frac{p}{p-1}$.
\fff
If one of the above conditions holds, then $\range{A_\ga}\subset\lpi$ and 
\iii 
\item[(i)] $A_\ga$ is a bounded linear operator from $\hh$ into $\lpi$;
\item[(ii)] its adjoint $A_\ga^*:\lqu\to\hh$ is given by 
\beeq{aggiunto}{A_\ga^*\phi=\int_X \phi(x)\ga_x\, \de\mu(x),}
where $\phi\in\lqu$  and the integral has to be understood in the weak
sense (if $p=\infty$ and $q=1$,  in Eq.\eqn{aggiunto} $A_\ga^*$ is
the restriction of the adjoint to $L^1(X,\mu)\subset L^{\infty}(X,\mu)^*$); 
\item[(iii)] $\Ga$ is a Carleman kernel and $A_\ga A_\ga^*=L_\Ga$, where
  $L_\Ga$ is the integral operator of kernel $\Ga$ given by Eq.\eqn{integrale}.
\fff
\end{proposition}
\bg{proof}~\\
\iii
\item[$1)\Rightarrow 2)$] To show that $A_\ga$ is a bounded
operator from $\hh$ into $\lpi$ we follow the proof of \citet{hiph}. 
Since $\ga$ is weakly $p$-integrable, $A_\ga$ is a linear
operator from $\hh$ to $\lpi$. We claim that it is closed. Indeed, let
$(v_n)_{n\in\nat}$ be a sequence that converges to
$v\in\hh$ and the sequence $(A_\ga v_n)_{n\in\nat}$ converges to
$\phi\in\lpi$. By construction, for $\mu$-almost all $x\in X$,
$$\lim_{n\to+\infty}(A_\ga v_n)(x)= \lim_{n\to+\infty}\scalh{v_n}{\ga_x}= \scalh{v}{\ga_x}=(A_\ga
v)(x).$$
By uniqueness of the limit, $A_\ga v=\phi$, so that $A_\ga$ is
closed. The closed graph theorem implies that $A_\ga$ is bounded.\\
We show Eq.\eqn{aggiunto}. Given
$\phi\in\lqu$ and $v\in\hh$, by assumption $A_\ga v\in\lpi$, so that the function 
$\phi
\overline{A_\ga v}=\phi(\cdot) \scalh{\ga(\cdot)}{v}$ is in $L^1(X,\mu)$. Since
$v$ is arbitrary, it follows that the function $\phi(\cdot)\, \ga(\cdot)$ is weakly
integrable and
$$\int_X \phi(x)\, \scalh{\ga_x}{v}\, \de\mu(x)=\scal{\phi}{A_\ga v}
=\scal{A_\ga^*\phi}{v}_\hh,$$
so that Eq.\eqn{aggiunto} holds. \\
The fact that $\Ga$ is a Carleman kernel follows observing that, for
all $x\in X$, $\Ga(\cdot,x)=A_\ga(\ga_x)(\cdot)$, which is in $\lpi$
by assumption. Moreover, given $\phi\in\lqu$,
Eqs.\eqn{Gamma} and\eqn{aggiunto} imply that
$$\int_X \Ga(x,t)\phi(t) \, \de\mu(t) = \scalh{A^*_\ga\phi}{\ga_x} =
(A_\ga A_\ga^*\phi)(x)\room \text{for}\ \maaa \ x\in X.$$
Hence, $\int_X \Ga(\cdot,t)\phi(t) \, \de\mu(t)$ is in $\lpi$, that is,
$\Ga$ is a Carleman $(q,p)$-bounded kernel, and $A_\ga
A_\ga^*=L_\Ga$.
\item[$2)\Rightarrow 1)$] Since $\Ga$ is measurable and $\hhg$
is separable, Proposition~\ref{misurabilita} ensures  the function 
 $\scalh{v}{\ga(\cdot)}$ be measurable for all $v\in\hh$. So it is
 enough to prove that  $|\scalh{v}{\ga(\cdot)}| \in
 L^p (X,\mu)$\\
To this aim, let $(A_n)_{n\in\nat}$ be the sequence of sets defined in
the proof of Proposition~\ref{op-misur}. For all $n\in\nat$ let
$D_n=\bigcup_{k=1}^n A_k$, which is of finite measure and $\ga$ is
bounded by $n$ on $D_n$. \\
Let now $\phi\in\lqu$, then the map $x\mapsto
\chi_{D_n}(x)\phi(x)\ga_x$ is clearly measurable and bounded by
$\chi_{D_n}n|\phi|\in L^1(X,\mu)$ since $D_n$ has finite measure.\\ 
Hence the linear operator $B_n:\lqu\to\hh$ 
$$ B_n\phi=\int_X \chi_{D_n}(x) \phi(x)\ga_x\,\de\mu(x)\room \phi\in\lqu,$$
where the integral is in the strong sense,  is well
defined and  bounded by  
$$ \norh{B_n\phi} \leq n \int_{D_n} |\phi(x)|\,\de\mu(x)
\leq n \mu (D_n)^{1/p} \norq{\phi}. $$
Since $\Ga$ is a bounded kernel, $L_\Ga$ is a bounded operator and the
following inequality holds
\begin{eqnarray*}
\norh{B_n\phi}^2 
& = & \int_X \chi_{D_n}(y) \overline{\phi(y)}
   \left(\int_X \chi_{D_n}(x) \phi(x) \scal{\ga_x}{\ga_y}_\hh \,\de\mu(x)\right)
   \,\de\mu(y) \\ 
& = & \scal{L_\Ga (\chi_{D_n}\phi)}{(\chi_{D_n}\phi)} \leq \noru{L_\Ga} \norq{\phi}^2,
\end{eqnarray*}
for all $\phi\in\lqu$, so that 
$$\sup_{n\in\nat}\noru{B_n} \leq M,$$
with $M=\sqrt{\noru{L_K}}$
and, hence, $\sup_n \noru{B_n^*} \leq M$. \\
Let now $v\in\hh$ and $\phi \in \lqu$, then
$$ \scal{\phi}{B_n^*v}= \scalh{B_n \phi}{v}=\int_X
\phi(x) \scalh{\chi_{D_n}(x)\ga_x}{v} \,\de\mu(x)\room \forall\phi\in
\lqu,$$
so $B_n^* v (x) = \chi_{D_n}(x)\scalh{v}{\ga_x}$ for $\mu$-almost all
$x\in X$ (in particular, $B_n^* v \in L^1 (X,\mu)$ if $q=\infty$), and 
\beeq{mio}{ \norp{B_n^* v} \leq \noru{B_n^*} \norh{v} \leq M \norh{v}.}
Since $\bigcup_{n\in\nat}D_n=X$, $\lim_{n\to+\infty}B_n^* v (x)=\scalh{v}{\ga_x}$
$\mu$-almost everywhere, and Eq.\eqn{mio} immediately ensures that 
$\scalh{v}{\ga(\cdot)} \in
L^\infty (X,\mu)$ for $p=\infty$. If $1\leq p <\infty$, the monotone convergence 
theorem
gives that the  map $x\mapsto |\scalh{v}{\ga_x}|^p$ is in $\luno$. 
\fff
\end{proof}
If  $\ga$ is a weakly measurable function such that 
$$\int_{X\times X} |\Ga(x,t)|^p\ \de(\mu\otimes\mu)(x,t)<+\infty,$$
Fubini theorem and Holder inequality imply that $\Ga$ is a
Carleman $(q,p)$-bounded kernel and, hence, $\ga$ is a weakly $p$-integrable 
function. \\
The following proposition discusses the problem of compactness, compare with
Corollary~\ref{compatto}. For finite measures the first statement is
due to \citet{pettis}. The second statement is well known. See, for example,
\citet{hasu} for a complete discussion about the compactness of
integral operators in $\ldue$.
\bg{proposition}
With the notation of Proposition~\ref{p-integrabilita}, 
if $\hhg$ is separable and $\ga$ is weakly integrable, then $A_\ga$ is
a compact operator in $\luno$.\\
If $p<\infty$ and $\ga$ is strongly $p$-integrable, then $A_\ga$ is
a compact operator in $\lpi$.
\end{proposition}
\bg{proof}
Assuming that $\hhg$ is separable and $\ga$ is weakly integrable, we prove that $A_\ga$ is
compact. First, suppose that $\ga$ takes only a countable number of values, that is,
$$\ga = \sum_{n\in I} \chi_{E_n} v_n,$$
where $I$ is denumerable, the sequence $(E_n)_{n\in I}$ is disjoint, 
each $E_n$ is of finite measure and $v_n\in\hh$. 
The condition that $\ga$ is weakly integrable implies that
$$\nor{A_\ga v}_{\lunot}=\sum_n \mu(E_n) |\scal{v}{v_n}| =
\nor{Tv}_{\ell_1},$$
where $T:\hh\to\ell_1$, $(Tv)_n=\mu(E_n) \scal{v}{v_n}$. Since $T$ is
a bounded operator from the Hilbert space $\hh$ to $\ell_1$, it is known
that $T$ is compact (see \citet{conway}). Hence $A_\ga$ maps weakly
convergent sequences into strongly convergent ones, so that it is compact.\\
Assume now that $\ga$ is arbitrary, we claim there exists $\ga_1$ and
$\ga_2$ such that $\ga_1$ takes
only a countable number of values, $\ga_2$ is strongly
integrable and $\ga=\ga_1+\ga_2$. To this aim, let $(A_n)_{n\in\nat}$ be the sequence 
of sets defined in
the proof of Proposition \ref{op-misur}. Given $n\in\nat$, since $A_n$ has finite
measure and $\ga$ is bounded by $n$ on it, $\chi_{A_n}\ga$ is strongly
integrable, so there is a map $\eta_n:A_n\to\hh$, which takes only a
finite number of values and
$$\int_{A_n}\norh{\ga(x)-\eta_n(x)}\ \de\mu(x) \leq \frac{1}{2^n}.$$
Since the sequence $(A_n)_{n\in\nat}$ is disjoint and measurable, the
map 
$$\ga_1=\sum_{n\in\nat} \chi_{A_n}\eta_n$$  
is well defined, measurable and takes only a countable number of
values. Let $\ga_2=\ga-\ga_1$, then
$$\int_X\norh{\ga_2(x)}\ \de\mu(x)=\sum_n \int_{A_n}\norh{\ga(x)-\eta_n(x)}\ \de\mu(x) 
\leq \sum_n\frac{1}{2^n}=2,$$
so that $\ga_2$ is strongly integrable and $\ga=\ga_1+\ga_2$.\\
Finally, since $\ga$ and $\ga_2$ are weakly integrable, so is $\ga_1$ and, by the
previous result, $A_{\ga_1}$ is compact. Since $\ga_2$ is strongly
integrable, it is easy to check that $A_{\ga_2}$ is compact (see the
second part of this proposition). The thesis
follows since $A_\ga = A_{\ga_1}+ A_{\ga_2}$.\\
The second statement of the proposition is standard. Indeed, if $\ga$ is a strongly
$p$-integrable map, that is, $\ga$ is measurable and
$$\int_X \norh{\ga_x}^p\, \de\mu(x) = \int_X \Ga(x,x)^{p/2}\, \de\mu(x) < +\infty,$$ 
then $A_\ga :\hh \rightarrow \lpi$ is a compact operator. Indeed, if
$(v_n)_{n\in\nat}$ is a sequence in $\hh$ which converges weakly to
$0$, then $\scalh{v_n}{\ga_x} \rightarrow 0$ for all $x\in X$. Since
$|\scalh{v_n}{\ga_x}|\leq \norh{v_n} \norh{\ga_x}$ and $\sup_n
\norh{v_n} <+\infty$, it follows by dominated convergence theorem
that $A_\ga v_n \rightarrow 0$ in $\lpi$. This shows that $A_\ga$ maps
weakly convergent sequences into norm convergent sequences, so $A_\ga$
is compact. 
\end{proof}
Example~\ref{esempio} below shows that, if $p>1$ and $\ga$ is only weakly $p$-integrable, $A_\ga$ can be a
non-compact operator in $\lpi$. However, if $\mu (X) <+\infty$, since $\lpi\subset\luno$,
$A_\ga$ is a compact operator from $\hh$ into $\luno$ \citep{hasu}. \\   
If $p=\infty$ and $\ga$ is essentially bounded, $A_\ga$ is not necessarily compact.

For $p=2$ we can compute $A_\ga^*A_\ga$, which is known as frame
operator in the context of frame theory (see, for example,
\citet{young}). The following result can be found in \citet{hasu}. 
\bg{corollary} 
With the notation of Proposition~\ref{p-integrabilita},
assume that $\ga$ is weakly square-integrable, then
\beeq{frame}{A_\ga^*A_\ga=\int_X
\scalh{\cdot}{\ga_x}\ga_x\, \de\mu(x),} where the integral
converges in the weak operator topology.

In particular, $A_\ga$ is a Hilbert-Schmidt operator if and only if
$\ga$ is strongly square-integrable, that is,
\beeq{hs}{\int_X\Ga(x,x)\,\de\mu(x)<+\infty.} 
If this last condition holds, the integral in Equation$\eqn{aggiunto}$
converges in norm and the integral in Equation$\eqn{frame}$ converges
in trace norm.

Finally, $\ga$ is a (strongly) square-integrable function if and only if 
$L_\Gamma$ is a trace class operator from  $\ldue$ in $\ldue$. 
\end{corollary}
\bg{proof} Eq.\eqn{frame} follows from Eq.\eqn{aggiunto} and the
definition of $A_\ga$.
To prove Condition\eqn{hs}, let $(e_n)_{n\in\nat}$ be a Hilbert basis of
$\hhg$. Since $A^*_{\ga}A_{\ga}$ is a positive operator and
$|\scalh{\ga(\cdot)}{e_n}|^2$ is a positive function,
by monotone convergence theorem, one has that
\mfi
\tx{Tr}\,(A^*_\ga A_\ga)& = & \sum_n \int_X |\scalh{e_n}{\ga_x}|^2\,
\de\mu(x) \\
& = & \int_X  \sum_n |\scalh{e_n}{\ga_x}|^2\,
\de\mu(x) \\
& = &  \int_X \scalh{\ga_x}{\ga_x} \,\de\mu(x) \\
& = & \int_X \Gamma(x,x) \,\de\mu(x)
\mff
and the thesis follows.
Finally, we prove the statements about Eqs.\eqn{aggiunto}
and\eqn{frame}. Indeed, since $\hhg$ is separable, by Proposition~\ref{misurabilita}, the 
map $\ga$
is measurable, so that the maps $x\mapsto \phi(x)\ga_x$ and
$x\mapsto\scalh{\cdot}{\ga_x}\ga_x$ are measurable as maps
taking values in $\hhg$ and in the separable Banach space of trace class operators on
$\hhg$, respectively. Moreover, 
\mfi
\norh{\phi(x)\ga_x} &= &
|\phi(x)|\sqrt{\Ga(x,x)} \in
L^1(X,\mu) \\
\nor{\scalh{\cdot}{\ga_x}\ga_x}_{\mathrm{tr}} & = & \Gamma(x,x) \in L^1(X,\mu),
\mff
where $\nor{\cdot}_{\mathrm{tr}}$ is the trace norm. So the integrals
are well defined and the  claims concerning convergence
follow.\\
The above result can be restated in the following
way. The function $\ga$ is (strongly)
square-integrable, that is, Eq.\eqn{hs} holds, if and only if
 $A_\ga^* A_\ga $ is a trace class operator on $\hh$ and, by polar decomposition,
this latter condition is equivalent to the fact that $A_\ga A_\ga^*$ is
of trace class on $\ldue$. This proves the last statement of the corollary. 
\end{proof}
The following example shows that the weakly square-integrability is not
sufficient for compactness of $A_\ga$ and that strong square-integrability is not
necessary.
\bg{example}\label{esempio}
Let $\hh$ be an infinite dimensional separable Hilbert space and $\mu$
a measure on $X$ such that there exists a
disjoint sequence $(X_n)_{n\in\nat}$ of measurable sets satisfying
$$0<\mu(X_n)=a_n <+\infty \room n\in\nat.$$
Let $(e_n)_{n\in\nat}$ be a basis of $\hh$ and $(\sigma_n)_{n\in\nat}$
be a positive sequence of $\runo$. \\
Define the map $\ga$ from $X$ to
$\hh$ as
$$\ga_x= \frac{1}{\sqrt{a_n}}\sigma_n e_n\room\tx{if}\ x\in X_n.$$ 
Clearly $\ga$ is measurable and, for all $v\in\hh$,
$$\scalh{v}{\ga_x}=  \frac{1}{\sqrt{a_n}}\sigma_n \scalh{v}{e_n}\room\tx{if}\ x\in 
X_n.$$
It follows that
$$\int_X |\scalh{v}{\ga_x}|^2 \de\mu(x) = \sum_n \sigma_n^2  |\scalh{v}{e_n}|^2.$$
Since the sequence $(|\scalh{v}{e_n}|^2)_{n\in\nat}$ is in $\ell_1$,
$\ga$ is a weakly square-integrable function if and only if the sequence
$(\sigma_n)_{n\in\nat}$ is bounded. \\
Assume from now on that $(\sigma_n)_{n\in\nat}$ is bounded. It follows that,
given $v\in \hh$,
$$A_\ga v= \sum_n \sigma_n \scalh{v}{e_n}\phi_n$$
where 
$$\phi_n(x)=\left\{\begin{array}{lr}
\frac{1}{\sqrt{a_n}} & x\in{X_n} \\
0 & x\not\in X_n
\end{array}\right.
$$
and the series converges in $\ldue$.\\
Since $(\phi_n)_{n\in\nat}$ is an orthonormal sequence of $\ldue$, the
above equation shows that $(e_n,\sigma_n^2)$ is the spectral
decomposition of $A_\ga^*A_\ga$. 
It follows that $A_\ga$ is a compact operator if and only if 
$$\lim_{n\to+\infty}\sigma_n = 0.$$
Finally, $A_\ga$ is a Hilbert-Schmidt operator if and only if 
$$\sum_{n=1}^{+\infty} \sigma_n^2 <+\infty.$$
\end{example}
We end the section with a simple proposition showing the connection
between integral operators and weakly integrable functions
(see \citet{hasu} for a discussion).
\bg{proposition}\label{op-int} 
With the notation of Proposition~\ref{op-misur}, let $K:X\times Y\to\cuno$ be
a Carleman $(p,2)$-bounded kernel with $1\leq p\leq\infty$ and let $\ga:X\to\lduet$
$$\ga_x=\overline{K(x,\cdot)}$$
for $\mu$-almost all $x\in X$, then $\ga$ is a weakly $p$-integrable function.
\end{proposition}
\bg{proof}
 By separability of $\lduet$, the map $x\mapsto \ga_x$ is measurable.
By definition of Carleman $(p,2)$-bounded kernel, for
$\mu$-almost all $x\in X$ $K(x,\cdot)$ is square-integrable with
respect to $\nu$ and 
$$\scal{v}{\ga_x}_{\lduenot} =  \int_Y K(x,y) v(y)\,\de\nu(y) =  (L_K
v)(x)\room\forall v\in\lduet.$$
Since $K$ is a $(p,2)$-bounded kernel, the range of $L_K$ is in
$\lpi$, so that $\gamma$ is a weakly $p$-integrable function.
\end{proof}

\section{Continuity}\label{sec_cont}
In this section we study the weak continuity of $\ga$. 
The following result is due to \citet[Proposition~24]{schwartz}, but our
proof is based on elementary tools.
\bg{proposition}\label{continuita}
Let $X$ be a locally compact space and $\hh$ a Hilbert space.  The following facts are 
equivalent:
\iii
\item  the map $\ga$ is weakly continuous;
\item the function $\Gamma$ is locally bounded and separately continuous.
\fff
If one of the above conditions holds, $\range{A_\ga}\subset\cc{X}$, the operator 
$A_\ga$ is
continuous from $\hh$ into $\cc{X}$, and
\beeq{genera1}{\Ker{A_\ga}=\overline{\rm span} \, \{ \ga_x\in \hh \vert x\in X 
\}^\perp = \hhg^\perp. } 
If $X$ is separable, then $\hhg$ is separable.
\end{proposition}
\bg{proof}{~}
\iii
\item[$1)\Rightarrow 2)$] By assumption, $A_\ga v$ is a continuous
function for all $v\in\hh$. 
Let $x\in X$, since $\Ga(\cdot,x)=A_\ga(\ga_x)(\cdot)$, clearly $\Ga$ is separately 
continuous. 
Fix now a compact set $C$. Since $A_\ga v=\scalh{v}{\ga(\cdot)}$ is
continuous, it is bounded on $C$ for all $v\in\hh$, so the
Banach-Steinhaus theorem ensures 
$$\norh{\ga_x}\leq M\room \forall x\in C,$$
for some constant $M>0$. Moreover, by Cauchy-Schwartz inequality,
$$|\Gamma(x,t)|=|\scalh{\ga_t}{\ga_x}|\leq \norh{\ga_x}\norh{\ga_t}\leq
M^2,$$
for all $t,x\in C$, so that $\Gamma$ is locally bounded.
\item[$2)\Rightarrow 1)$] Let $w=\sum_{i=1}^n a_i \ga_{x_i}$.
The fact that $\Ga$ is separately continuous and Eq.\eqn{Gamma}
imply that the function $A_\ga w=\sum_{i=1}^n a_i \Ga(\cdot,x_i)$
is continuous. Let now $v\in\hhg$ and $x_0\in X$, we prove that
$A_\ga v$ is continuous at $x_0$. Let $C$ be a
compact neighborhood $C$ of $x_0$ and 
$$M=\sup_{x\in C}\norh{\ga_x}=\sup_{x\in C}\sqrt{\Ga(x,x)},$$
where $M$ is finite due to locally boundedness of $\Ga$. Fixed
$\eps>0$, there is a finite linear combination $w=\sum a_i \ga_{x_i}$
such that $\norh{v-w}\leq\eps$. By the above observation,
$A_\ga w$ is continuous, so, possibly replacing $C$
with a smaller neighborhood,
$|(A_\ga w)(x)-(A_\ga w)(x_0)|\leq\eps $ for all $x\in C$. Then 
\mfi
|(A_\ga v)(x)-(A_\ga v)(x_0)| & \leq &
|(A_\ga w)(x)-(A_\ga w)(x_0)|+ \norh{\ga_x-\ga_{x_0}}\norh{w-v} \\
& \leq &\eps (1+ 2M)\room\forall x\in C.
\mff
Finally, if  $v\in\hhg^\perp$, $A_\ga v=0$, so that $\ga$ is weakly continuous.
\fff
If one of the above equivalent conditions holds, $A_\ga$ is a
bounded operator from $\hh$ into $\cc{X}$. Indeed, let $C$ be a compact set, by 
locally boundedness of
$\Ga$ and $\nor{\ga_x}^2=\Ga(x,x)$, there is a constant $M>0$ such that
$\norh{\ga_x}\leq M$ for all $x\in C$. Finally,
$$
\sup_{x\in C} |(A_\ga v)(x)| =  \sup_{x\in C} |\scalh{v}{\ga_x}| 
\leq  \norh{v} \sup_{x\in C}\norh{\ga_x} \leq M \norh{v},
$$
so that $A_\ga$ is bounded. The fact that $\Ker{A_\ga}=\hhg^\perp$ is  clear.
Finally, assume that $X$ is separable, we prove that $\hhg$ 
is separable. Indeed, let $X_0$ be a dense denumerable subset of $X$ and define
$$\hh_0=\overline{\tx{span}}\set{\ga_x\,|\,x\in X_0}\subset \hhg.$$
Clearly, $\hh_0$ is separable. We claim that $\hh_0=\hhg$.
Choose $v\in\hh_0^\perp$ so that  $\left(A_\gamma v\right)(x) = \scalh{v}{\ga_x}=0$ 
for all
$x\in X_0$. Since $A_\gamma v$ is continuous and $X_0$ is dense,
$\scalh{v}{\gamma_x} =0$ for all $x\in X$ so that $v\in\hhg^\perp$. It
follows that $\hh_0^\perp\subset\hhg^\perp$ so that $\hhg\subset\hh_0$ and the claim 
follows.
\end{proof}
The following corollary characterizes the strong continuity of
$\ga$ (see Proposition~24 of~\citet{schwartz}).
\bg{corollary}\label{compatto}
With the assumptions of Proposition~\ref{continuita}, the following facts are 
equivalent:
\iii
\item the function $\Gamma$ is continuous;
\item the function $\Gamma$ is continuous on the diagonal of $X\times X$;
\item the map $\ga$ is continuous;
\item the map $\ga$ is weakly continuous and the restriction of $\Ga$
  on the diagonal $x\mapsto\Ga(x,x)$ is continuous;
\item the operator $A_\ga:\hh\to\cc{X}$ is compact.
\fff
\end{corollary}
\bg{proof}{~}
\iii
\item[$1)\Rightarrow 2)$] Trivial.
\item[$2)\Rightarrow 3)$] It follows observing that Eq.\eqn{Gamma}
  gives 
$$ \norh{\ga_x-\ga_t}^2=\Ga(x,x)+\Ga(t,t)- \Ga(x,t) - \Ga(t,x)\room x,t\in X.$$
\item[$3)\Rightarrow 1)$] Use the fact that the scalar product is
  continuous and Eq.\eqn{Gamma}.
\item[$3)\iff 4)$] It is a restatement
of the equivalence between strong convergence and weak convergence
plus convergence  of norms (so called $\hh$-property).
\item[$3)\iff 5)$] Let
$B_1$ be the closed unit ball in $\hh$. The operator  $A_\ga$ is compact if and
only if $A_\ga(B_1)$ is compact in $\cc{X}$. Since a locally compact
Hausdorff space is regular, by Ascoli theorem \citep{kelley} the
set of functions $A_\ga(B_1)$ is compact if and 
only if $A_\ga(B_1)$ is closed, the set $(A_\ga(B_1))(x)$ is bounded
for all $x\in X$, and $A_\ga(B_1)$ is equicontinuous. First two
conditions are always satisfied since $A_\ga$ is continuous from the
Hilbert space $\hh$ to $\cc{X}$. Moreover, given $x,t\in X$, we have
that 
\mfni
\norh{\ga_x-\ga_t} & = & \sup_{v\in B_1}|\scalh{v}{\ga_x-\ga_t}| \nonumber \\
& = & \sup_{v\in B_1}|(A_\ga v)(x)-(A_\ga v)(t)| \nonumber\\
& = & \sup_{f\in A_\ga(B_1)}|f(x)-f(t)|. \label{equi} 
\mfnf
Eq.\eqn{equi} shows that $A_\ga(B_1)$ is equicontinuous if and
only if map $\ga$ is continuous.
\fff  
\end{proof}
Let  now $X$ be a locally compact
second countable Hausdorff space endowed with a positive Radon measure
$\mu$. Assume that $\ga$ is a weakly continuous map from $X$ into a Hilbert
space $\hh$ and regard $A_\ga$ as an operator from $\hh$ into $\emme$. 
Since $\hhg$ is separable,  
Proposition~\ref{misurabilita} holds, so that $
\ker{A_\ga}=\hhm^\perp$. On the other hand, since  $f\in\cc{X}$ is equal to $0$ in 
$\emme$ if and only if
$f(x)=0$ for all $x\in\supp{\mu}$, 
$$ \Ker A_\ga = \{ v \in \hh \vert \scalh{v}{\ga_x} = 0 \,\forall x \in \supp{\mu} \}
= \overline{\rm span}\,\{ \ga_x \vert x \in \supp{\mu} \}^\perp , $$
hence
\beeq{kernel-bello}{\hhm=\overline{\rm span}\, \{\ga_x\,\vert\, x\in\supp{\mu} \}. } 
Finally, assume $X$ compact, so that $\mu$ is finite. Since
$\ga$ is weakly continuous, it is bounded and, hence, strongly
$p$-integrable, so the operator $A_\ga$ is always compact as a map in
$L^p(X,\mu)$. However, in order $A_\ga$ be compact as a map in $\cc{X}$, it is
necessary (and sufficient) that $\ga$ is strongly continuous.

\section{Mercer theorem}\label{sec_mercer}
In this section we characterize those reproducing kernel Hilbert
spaces that are subspaces of $\ldue$ in terms of the spectral
decomposition of the integral operator $L_\Ga$. We recall the
definition of complex Radon measure \citep{dieu}. Let $R>0$ and  $\ccc{(0,R]}$ the
space of compactly supported functions on the interval $(0,R]$, a {\em complex Radon 
measure}
on $(0,R]$ is a complex linear form on $\ccc{(0,R]}$ such that its
restriction to $\ccc{[a,R]}$ is bounded for all $0<a<R$. 
If $\rho$ is a complex Radon measure, there is a unique positive Radon measure
$|\rho|$  and a complex measurable function $h$ on $(0,R]$,
such that
$$\rho(\phi)=\int_{(0,R]} \phi(\la)h(\la)\de|\rho|(\la)\room \phi\in\ccc{(0,R]}.$$
and $|h(\la)|=1$ for all $\la\in (0,R]$.
If $\phi\in L^1((0,R],|\rho|)$, the integral with respect to $\rho$ is
defined by
$$ \int_{(0,R]} \phi(\la) \de\rho(\la) = \int_{(0,R]} \phi(\la) h(\la) 
\de|\rho|(\la).$$
The measure $|\rho|$ is called the {\em absolute value} of $\rho$.
\bg{proposition}\label{mercer}
Let $X$ be a locally compact second countable Hausdorff space endowed
with a positive Radon measure $\mu$ such that 
$\supp{\mu}=X$. Let $\hh$ be a reproducing kernel Hilbert space with
kernel $\Ga$. The following conditions are equivalent.
\iii
\item $\hh$ is a subspace of $\ldue\cap \cc{X}$;
\item the kernel $\Ga$ is locally bounded, separately continuous, and $(2,2)$-bounded.
\fff
If one of the above assumptions holds, the integral operator $L_\Ga$ of
kernel $\Ga$ is a positive operator with  spectral decomposition 
$$ L_\Ga = \int_{[0,R]} \lambda \de P(\lambda) \room
R=\nor{L_\Ga} $$
and
\mfni
\hh & = & \set{v\in\ldue\cap\cc{X} \,\vert\, \int_{[0,R]}\frac{1}{\la} 
  \scal{\de P(\la)v}{v}_{\ldueXnot} <+\infty} \label{mercer1} \\
\norh{v}^2 & = & \int_{[0,R]} \frac{1}{\lambda} 
\scal{\de P(\la)v}{v}_{\ldueXnot}\room \forall v\in\hh.\label{mercer2}
\mfnf 
If $x,y\in X$ there is a complex Radon measure
$\rho_{x,y}$ on $(0,\R]$ such that
\beeq{mercer3}{\Ga(x,y)= \int_{(0,R]} \la \de \rho_{x,y}(\la)}
and, for all $E\in\mathcal{B}((0,\R])$ such that
$\overline{E}\subset(0,\R]$, then
\beeq{mercer4} {\int_E \de\rho_{x,y}(\la)=\sum_{n\in I} \phi_n(x) 
\overline{\phi_n(y)}}
where $\phi_n\in\hh$ and $(\phi_n)_{n\in I}$ is an  orthonormal basis in $\ldue$ for 
the range of $P(E)$.\\
If $x\in X$, $\rho_{x,x}$ is a positive Radon measure on
$(0,\R]$ such that
\beeq{mercer5} {\rho_{x,x} ((0,\R])=\sum_{n\in I} |\phi_n(x)|^2, }
where $\phi_n\in\hh$ and $(\phi_n)_{n\in I}$ a basis in $\ldue$ for
  $\Ker{L_\ga}^\perp$. In particular $\rho_{x,x}$ is finite if and
  only if $\Ga(\cdot,x)\in\range{L_\Ga}$.
\end{proposition}
\bg{proof}
Proposition~\ref{range} ensures that the inclusion of $\hh$
into $\cuno^X$ is  the operator $A_\ga$ associated with the map $\ga:X\to\hh$
$$\ga_x=\Ga(\cdot,x)\room x\in X.$$
Proposition~\ref{continuita} imply that
$\hh\subset\cc{X}$ if and only if  $\Ga$ is separately continuous and
locally bounded. In particular, since $X$ is separable, both the above
conditions and Eq.\eqn{genera} ensure that $\hh$ is separable. So
Proposition~\ref{p-integrabilita} with $p=2$ gives that
$\hh\subset\ldue$ if and only if $\Ga$ is a $(2,2)$-bounded kernel.
Hence the equivalence of Condition~1 and Condition~2 is now clear.
Assume one of them and regard the inclusion $A_\ga$ as an operator from
$\hh$ into $\ldue$. By Eqs.\eqn{kernel-bello},\eqn{kernel},
$A_\ga$ is injective and the polar decomposition of the adjoint $A_\ga^*$ gives
\beeq{III1}{A_\ga^*=W(A_\ga A_\ga^*)^{\frac 12}}
where $W$ is a partial isometry from $\ldue$ to $\hh$ with kernel being equal to the 
kernel of $A_\ga^*$
and with range being equal to the closure of the range of $A_\ga^*$.

Proposition~\ref{p-integrabilita} states that $A_\ga A_\ga^*=L_\Ga$, so
that $L_\Ga$ is a positive operator,
the kernel of $W$ is the kernel of $L_\Ga$ and, since $A_\ga$ is
injective, $W$ is surjective, that is,  $WW^*$ is the
identity. Equation\eqn{III1} gives
\beeq{III2}{A_\ga =L_\Ga^{\frac 12} W^*.}
Since $A_\ga$ is the inclusion of $\hh$ into
$\ldue\cap\cc{X}$  and is injective, Eq.\eqn{III2} implies that
 we can identify $\hh$ with $\range{L_\Ga^{\frac 12}}$ and, by means of spectral 
theorem,
Eq.\eqn{mercer1} follows. In particular, for all $v\in\hh$, 
$v\in(\Ker{L_\Ga})^\perp=(\range{P(\{0\})})^\perp$.

Let now $v\in\hh$, then $v=L_\Ga^{\frac 12} W^*v$ and
\mfi
 \int_{[0,R]}\frac{1}{\la} 
  \scal{\de P(\la) v}{ v}_{\ldueXnot} & = & 
 \int_{[0,R]}\frac{1}{\la}  \scal{\de P(\la)L_\Ga^{\frac 12} W^*v }{
   L_\Ga^{\frac 12} W^*v}_{\ldueXnot} \\
 & = & \int_{[0,R]} \scal{\de P(\la) W^*v }{W^*v}_{\ldueXnot} \\
& = &  \scal{ W^*v }{W^*v}_{\ldueXnot} \\
& = &  \scalh{v }{v}, 
\mff
where we used that $P(\{0\})W^*=0$. So 
Eq.\eqn{mercer2} follows.

Given $x,y\in X$, let $\rho_{x,y}$ be the linear form on $\ccc{(0,R]}$
given by
$$ \rho_{x,y}(\phi)= \int_{(0,R]} \frac{1}{\la} \phi(\la) \scal{\de
  P(\la)W^*\ga_y}{W^*\ga_x}_{\ldueXnot}.$$
Since $\scal{\de P(\la)W^*\ga_y}{W^*\ga_x}$ is a bounded complex measure and
$\phi$ has compact support, $\rho_{x,y}$ is well defined and is a
complex Radon measure. By definition, $|\rho_{x,y}|$ has density $\frac{1}{\la}$ with 
respect to the
positive bounded measure $|\scal{\de P(\la)W^*\ga_y}{W^*\ga_x}|$. Since the function 
$1$
is integrable with respect to $|\scal{\de P(\la)W^*\ga_y}{W^*\ga_x}|$,
the function $\la$ is integrable with respect to
$|\rho_{x,y}|$ and, hence, with respect to $\rho_{x,y}$. In particular,
\mfi
 \int_{(0,R]} \la \de \rho_{x,y}  & =  &  \int_{(0,R]} \scal{\de
   P(\la)W^*\ga_y}{W^*\ga_x}_{\ldueXnot}  \\
& = & \int_{[0,R]} \scal{\de
   P(\la)W^*\ga_y}{W^*\ga_x}_{\ldueXnot} \\
& = &  \scal{W^*\ga_y}{W^*\ga_x}_{\ldueXnot} =   \scalh{\ga_y}{\ga_x},
\mff
so that Eq.\eqn{mercer3} follows by definition of $\Ga$. 

Let now $E\in\mathcal{B}((0,\R])$ be such that
$\overline{E}\subset(0,\R]$. This last fact and the spectral theorem
imply that $\range{P(E)}\subset  \range{L_\Ga^{\frac 12}}=\hh$. Hence,
there is a sequence $(\phi_n)_{n\in I}$ of $\hh$ such that
$(\phi_n)_{n\in I}$ is an  orthonormal basis for $\range{P(E)}$.
Since $\overline{E}\subset(0,\R]$, $\chi_E$ is integrable with respect to $\rho_{x,y}$ 
and
\mfni
 \int_{(0,R]} \chi_E(\la) \de \rho_{x,y} & = &  \int_{[0,R]} 
\frac{\chi_E(\la)}{\la}\scal{\de
   P(\la)W^*\ga_y}{W^*\ga_x}_{\ldueXnot} \nonumber\\
 & = &  \int_{[0,R]} \frac{1}{\la}
 \scal{\de P(\la)P(E)W^*\ga_y}{P(E)W^*\ga_x}_{\ldueXnot}
 \nonumber\\
& = &  \scal{L_\Ga^{-\frac 12}P(E)W^*\ga_y}{L_\Ga^{-\frac 
12}P(E)W^*\ga_x}_{\ldueXnot},\label{V1}
\mfnf
where $P(E)W^*\ga_x$ and $P(E)W^*\ga_y$ are in $\range{L_\Ga^{\frac 12}}$.

Let now $J$ a finite subset of $I$, taking into account the properties
of the sequence $(\phi_n)_{n\in I}$,{\setlength\arraycolsep{0pt}
\begin{eqnarray*}
&&\sum_{n\in J} \scal{L_\Ga^{-\frac 12}P(E)W^*\ga_y}{\phi_n}_{\ldueXnot}
\scal{\phi_n}{L_\Ga^{-\frac 12}P(E)W^*\ga_x}_{\ldueXnot}  \\
&& \qquad \qquad  =  \sum_{n\in J} \scal{W^*\ga_y}{L_\Ga^{-\frac 12}\phi_n}_{\ldueXnot} 
\scal{L_\Ga^{-\frac 12}\phi_n}{W^*\ga_x}_{\ldueXnot} \\
&& \qquad \qquad (\phi_n=L_\Ga^{\frac 12}W^*\phi_n) \\
&& \qquad \qquad = \sum_{n\in J} \scal{W^*\ga_y}{ W^*\phi_n}_{\ldueXnot} 
\scal{ W^*\phi_n}{W^*\ga_x}_{\ldueXnot} \\
&& \qquad \qquad  =  \sum_{n\in J} \scalh{\ga_y}{\phi_n}\scalh{\phi_n}{\ga_x} \\
&& \qquad \qquad (Eq.\eqn{valuta}) \\
&& \qquad \qquad =  \sum_{n\in J} \phi_n(x) \overline{\phi_n(y)}.
\end{eqnarray*}}
Eq.\eqn{mercer4}
follows observing that the series
$$ \sum_{n\in I} \scal{L_\Ga^{-\frac 12}P(E)W^*\ga_y}{\phi_n}_{\ldueXnot} 
\scal{\phi_n}{L_\Ga^{-\frac 12}P(E)W^*\ga_x}_{\ldueXnot}$$
is summable with sum $\scal{L_\Ga^{-\frac 12}P(E)W^*\ga_y}{L_\Ga^{-\frac
    12}P(E)W^*\ga_x}_{\ldueXnot}= \rho_{x,y}(E)$, by means of Eq.\eqn{V1}.
 
Finally, if $x\in X$, clearly $\rho_{x,x}$ is a positive Radon measure
on $(0,R]$ having density $\frac{1}{\la}$ with respect to the positive bounded
measure $\scal{\de P(\la)W^*\ga_x}{W^*\ga_x}$. Hence $\rho_{x,x}$ is
bounded if and only if
$W^*\ga_x\in\range{L_\Ga^{\frac 12}}$. Equation~\eqn{III2} implies that
this condition is equivalent to $\ga_x\in\range{L_\Ga}$ and, if it is
satisfied,
\beeq{III4}{\rho_{x,x}((0,R]) = \nor{L_\Ga^{-\frac 12}W^*\ga_x}^2_{\ldueXnot}.}  
Let $(\phi_n)_{n\in\nat}$ be a sequence in $\hh$ such that
$(\phi_n)_{n\in\nat}$ is a basis in $\ldue$ of $\Ker{L_\Gamma}^\perp$.
Given $N\in\nat$, reasoning as above, 
$$ \sum_{n=1}^N |\phi_n(x)|^2  =  \sum_{n=1}^N 
\left|\scal{L_\Ga^{-\frac 12 }\phi_n}{W^*\ga_x}_{\ldueXnot}\right|^2 .$$
The series in the right side converges if and only if
$W^*\ga_x\in\range{L_\Ga^{\frac 12}}$ and, if it is so, its sum is
$\nor{L_\Ga^{-\frac 12}W^*\ga_x}^2_{\ldueXnot}$. Eq.\eqn{III4} implies Eq.\eqn{mercer5}.
\end{proof} 
Equation\eqn{mercer1} allows us to identify the elements of $\hh$ with
the only continuous functions on $X$ whose equivalence class belongs to
the range of $L_\Ga^{\frac{1}{2}}$, extending a result of
\citet{cusm}. The assumptions that $\supp{\mu}=X$ and $\hh\subset\cc{X}$ ensure that
the identification between functions in $\hh$ and equivalence classes
in $\ldue$ is well defined. With this identification,
Eq.\eqn{mercer2} implies that $L_\Ga^{\frac 12}$ is a unitary
operator from $\Ker{L_\Ga}^\perp$ onto $\hh$. \\
If $\supp{\mu}\neq X$, the statements of the above proposition hold
replacing $X$ with $\supp{\mu}$ and $\hh$ with $\overline{\rm span}\,\{ \ga_x
\vert x\in \supp{\mu} \}$. We use the assumption that $\ga$ is
weakly continuous only in two steps: to show the separability of
$\hh$ and to identify elements of $\hh$, which are
functions on $X$, with elements in $\ldue$, which are equivalence
classes. 

If the integral operator $L_\Ga$ has a pure point spectrum 
$$L_\Ga  =  \sum_{n\in I} \la_n \scal{\cdot}{\phi_n}\phi_n,$$
where $\la_n\geq 0$ and $(\phi_n)_{n\in\nat}$ is a basis of $\ldue$,
it is possible to choose $\phi_n\in\hh$ for all $\la_n>0$ and, with
this choice, Proposition~\ref{mercer} gives that
\mfi
\hh & = & \set{v\in\ldue\cap\cc{X}\,|\,\sum_n \frac{1}{\la_n}
  |\scal{v}{\phi_n}_{\ldueXnot}|^2 <+\infty}, \\
\nor{v}_\hh^2 & = & \sum_n \frac{1}{\la_n} |\scal{v}{\phi_n}_{\ldueXnot}|^2, \\
\Gamma(x,y) & = & \sum_n \la_n \phi_n(x)\overline{\phi_n(y)}.
\mff
The last series converges absolutely and, by Dini theorem and
Condition~4 of Corollary~\ref{compatto}, converges
uniformly on compact subsets if and only if $\ga$ is strongly continuous.


\end{document}